\documentclass[11pt]{article}
\usepackage{amsfonts,amsmath,amsxtra}
\usepackage{latexsym}
\usepackage{amssymb}

\def\hybrid{\topmargin 0pt      \oddsidemargin 0pt
        \headheight 0pt \headsep 0pt
        \textwidth 16.5cm
        \textheight 24cm
        \marginparwidth 0.0in
        \parskip 5pt plus 1pt   \jot = 1.5ex}
\catcode`\@=11
\def\marginnote#1{}
\newcount\hour
\newcount\minute
\newtoks\amorpm
\hour=\time\divide\hour by60 \minute=\time{\multiply\hour by60
\global\advance\minute by-\hour}
\edef\standardtime{{\ifnum\hour<12 \global\amorpm={am}%
        \else\global\amorpm={pm}\advance\hour by-12 \fi
        \ifnum\hour=0 \hour=12 \fi
      \number\hour:\ifnum\minute<10 0\fi\number\minute\the\amorpm}}
\edef\militarytime{\number\hour:\ifnum\minute<10 0\fi\number\minute}

\def\draftlabel#1{{\@bsphack\if@filesw {\let\thepage\relax
   \xdef\@gtempa{\write\@auxout{\string
      \newlabel{#1}{{\@currentlabel}{\thepage}}}}}\@gtempa
   \if@nobreak \ifvmode\nobreak\fi\fi\fi\@esphack}
        \gdef\@eqnlabel{#1}}
\def\@eqnlabel{}
\def\@vacuum{}
\def\draftmarginnote#1{\marginpar{\raggedright\scriptsize\tt#1}}

\def\draft{\oddsidemargin -0.1truein
        \def\@oddfoot{\sl preliminary draft \hfil
        \rm\thepage\hfil\sl\today\quad\militarytime}
        \let\@evenfoot\@oddfoot \overfullrule 3pt
        \let\label=\draftlabel
        \let\marginnote=\draftmarginnote
\def\@eqnnum{{\rm (\theequation)}
\rlap{\kern\marginparsep\tt\@eqnlabel}%
\global\let\@eqnlabel\@vacuum}  }

\newfont{\Bbbb}{msbm7 scaled 1\@ptsize00}
\newcommand{\zs}{\raise-1pt\hbox{$\mbox{\Bbbb Z}$}}

scaled 1\@ptsize00

\font\sevenmsa=msam6 
\newfam\msafam
\textfont\msafam=\sevenmsa
\def\hexnumber@#1{\ifnum#1<10 \number#1\else
\ifnum#1=10 A\else\ifnum#1=11 B\else\ifnum#1=12 C\else \ifnum#1=13
D\else\ifnum#1=14 E\else\ifnum#1=15 F\fi\fi\fi\fi\fi\fi\fi}
\def\msa@{\hexnumber@\msafam}
\def\llcorner{\delimiter"4\msa@78\msa@78 }
\def\lrcorner{\delimiter"5\msa@79\msa@79 }
\mathchardef\blacktriangleright="3\msa@49
\mathchardef\blacktriangleleft="3\msa@4A \font\tenmsb=msbm10 scaled
1\@ptsize00
\newfam\msbfam
\textfont\msbfam=\tenmsb \scriptfont\msbfam=\tenmsb

\newdimen\Squaresize \Squaresize=14pt
\newdimen\Thickness \Thickness=0.5pt

\def\Square#1{\hbox{\vrule width \Thickness
   \vbox to \Squaresize{\hrule height \Thickness\vss
      \hbox to \Squaresize{\hss#1\hss}
   \vss\hrule height\Thickness}
\unskip\vrule width \Thickness} \kern-\Thickness}

\def\Vsquare#1{\vbox{\Square{$#1$}}\kern-\Thickness}

\def\numberbysection{\@addtoreset{equation}{section}
        \def\theequation{\thesection.\arabic{equation}}}
\numberbysection

\renewcommand{\theequation}{\thesection.\arabic{equation}}
\def\titlepage{\@restonecolfalse\if@twocolumn\@restonecoltrue\onecolumn
     \else \newpage \fi \thispagestyle{empty}\c@page\z@
        \def\thefootnote{\fnsymbol{footnote}} }

\def\endtitlepage{\if@restonecol\twocolumn \else  \fi
        \def\thefootnote{\arabic{footnote}}
        \setcounter{footnote}{0}}  
\relax

\hybrid
\makeatletter
\newdimen\normalarrayskip            
\newdimen\minarrayskip               
\normalarrayskip\baselineskip \minarrayskip\jot
\newif\ifold             \oldtrue            \def\new{\oldfalse}
\def\arraymode{\ifold\relax\else\displaystyle\fi}
\def\eqnumphantom{\phantom{(\theequation)}} 
\def\@arrayskip{\ifold\baselineskip\z@\lineskip\z@
     \else
     \baselineskip\minarrayskip\lineskip1\baselineskip\fi}

\def\@arrayclassz{\ifcase \@lastchclass \@acolampacol \or
\@ampacol \or \or \or \@addamp \or
   \@acolampacol \or \@firstampfalse \@acol \fi
\edef\@preamble{\@preamble
  \ifcase \@chnum
     \hfil$\relax\arraymode\@sharp$\hfil
     \or $\relax\arraymode\@sharp$\hfil
     \or \hfil$\relax\arraymode\@sharp$\fi}}

\def\@array[#1]#2{\setbox\@arstrutbox=\hbox{\vrule
     height\arraystretch \ht\strutbox
     depth\arraystretch \dp\strutbox
width\z@}\@mkpream{#2}\edef\@preamble{\halign \noexpand\@halignto
\bgroup \tabskip\z@ \@arstrut \@preamble \tabskip\z@ \cr}%
\let\@startpbox\@@startpbox \let\@endpbox\@@endpbox
    \if #1t\vtop \else \if#1b\vbox \else \vcenter \fi\fi
  \bgroup \let\par\relax
  \let\@sharp##\let\protect\relax
  \@arrayskip\@preamble}
\def\eqnarray{\stepcounter{equation}%
              \let\@currentlabel=\theequation
              \global\@eqnswtrue
              \global\@eqcnt\z@
              \tabskip\@centering              
              \let\\=\@eqncr
              $$%
            \halign to \displaywidth  \bgroup
             \eqnumphantom \@eqnsel
      \hskip\@centering                               
    $\displaystyle  \tabskip\z@ {##}$%
    &\global\@eqcnt\@ne \hskip 2\arraycolsep
         $ \displaystyle  \arraymode{##}$\hfil
    &\global\@eqcnt\tw@ \hskip 2\arraycolsep
         $\displaystyle\tabskip\z@{##}$\hfil
         \tabskip\@centering
    &{##}\tabskip\z@\cr}
\makeatother

\newcommand{\CC}{{\mathbb{C}}}

\def\IC{\mathbb{C}}

\def\IP{\mathbb{P}}

\def\IR{\mathbb{R}}

\def\CC {\mathcal{C}}

\def\CM {\mathcal{M}}

\def\CS {\mathcal{S}}

\def\CU {\mathcal{U}}
\def\CV {\mathcal{V}}
\def\CW {\mathcal{W}}

\def\frak{\mathfrak}

\def\we{\raise-1pt\hbox{$\,\stackrel{\wedge}{,}\,$}}
\def\tr{{\rm tr}\,}
\def\pr {\partial}

\def\ov {{\overline}}

\def\be {{\beta}}
\def\ga{{\gamma}}
\def\s {{\sigma}}

\def\la{\lambda}

\def\e{\epsilon}
\def\pr {\partial}

\newcommand{\ul}{{\underline{\la}}}
\newcommand{\ug}{{\underline{\gamma}}}

\newcommand{\ux}{{\underline{x}}}

\newcommand{\gl}{{\mathfrak{gl}}}

\newcommand{\Fl}{{\mathrm{ Fl}}}
\newcommand{\Gr}{{\mathrm{ Gr}}}
\newcommand{\h}{{\hbar}}

\newcommand{\Ind}{{\mathrm{ Ind}}}

\newtheorem{te}{Theorem}[section]
\newtheorem{de}{Definition}[section]
\newtheorem{prop}{Proposition}[section]

\newtheorem{lem}{Lemma}[section]

\newtheorem{rem}{Remark}[section]

\newcommand\bqa{\begin{eqnarray}}
\newcommand\eqa{\end{eqnarray}}
\def\be{\begin{eqnarray}\new\begin{array}{cc}}
\def\ee{\end{array}\end{eqnarray}}

\def\beq{\begin{equation}}
\def\eeq{\end{equation}}
\def\bse{\begin{subequations}}                
\def\ese{\end{subequations}}
\def\bp{\begin{pmatrix}}
\def\ep{\end{pmatrix}}

\def\h{\hbar}
\def\i{\imath}

\newcommand\rk{\operatorname{rank}}
\newcommand\opp{\operatorname{opp}}

\def\stack#1#2{\raise0.7pt\hbox{$\mathrel{\mathop{#2}\limits^{#1}}$}}
\def\tr{\triangleright}
\def\tl{\triangleleft}
\def\sem{\mathsurround=0pt \raise1pt
\hbox{$\scriptscriptstyle>\!\!$}\:\!\!\tl}
\def\mes{\mathsurround=0pt \tr\!\:\!\raise0.8pt
\hbox{$\scriptscriptstyle\!\!<$}\,}
\def\]{\mathsurround=0pt ]\raise-2pt\hbox{$_\ast$}}

\def\<{\langle}
\def\>{\rangle}
\def\ov{\overline}

\newcounter{pac}[section]

\newcounter{pacc}[subsection]
\newcommand{\npaa}{\addtocounter{pacc}{1} \noindent {\bf
\arabic{section}.\arabic{subsection}.\arabic{pacc}}\,\,\,}
  0

\begin{document}
\title{\bf On parabolic Whittaker functions}
\author{\sc\Large Sergey Oblezin\footnote{The work is
partially supported by P. Deligne's 2004 Balzan Prize in
Mathematics, and by the RFBR 09-01-93108-NCNIL-a Grant.}}

\date{}
\maketitle

\begin{abstract}\noindent We derive a Mellin-Barnes
integral representation for solution to generalized (parabolic)
quantum Toda lattice introduced in \cite{GLO}, which presumably
describes the $S^1\times U_N$-equivariant Gromov-Witten invariants
of Grassmann variety.
\end{abstract}

\section*{Introduction}

The $\gl_N$-Whittaker functions, being solutions to the quantum
cohomology D-module $QH^*(\Fl_N)$ of the complete flag variety
$\Fl_N=GL_N(\IC)/B$, describe the corresponding equivariant
Gromov-Witten invariants of $\Fl_N$ (see \cite{Giv1}, \cite{Giv2}
and references therein). However, the Givental's approach to
representation theory description of quantum cohomology of
homogeneous spaces is inapplicable to generic incomplete flag
variety $\Fl_{m_1,\ldots,m_k}$, since no relevant Whittaker model
(Toda lattice) associated with an incomplete flag variety was known.

From the other hand, in \cite{HV} it was conjectured a description
of quantum cohomology of Grassmannians in terms of (non-Abelian)
gauged topological theories, together with a period-type integral
representation for the corresponding generating function.

Recently, in \cite{GLO} a generalization of the $\gl_N$-Whittaker
function to the case of the Grassmann variety
$\Gr_{m,\,N}=GL_N(\IC)/P_m,\,1\leq m<N$ is proposed. Namely, in
\cite{GLO} it is defined a Toda-type D-module and its solution
$\Psi^{(m,N)}_{\la_1,\ldots,\la_N}(x_1,\ldots,x_N)$ (referred to as
$\Gr_{m,\,N}$-Whittaker function), such that after specialization
$x_2=\ldots=x_N=0$ the symbols of this D-module reproduce the small
quantum cohomology algebra $qH^*(\Gr_{m,\,N})$ according to
\cite{AS} and \cite{K}. Conjecturally, the constructed generalized
Whittaker function describe the equivariant Gromov-Witten invariants
of $\Gr_{m,N}$, and in \cite{GLO} this conjecture is verified in
particular case of projective space $\IP^{N-1}=\Gr_{1,N}$.

In this note we construct Mellin-Barnes type integral representation
of the specialized $\Gr_{m,N}$-Whittaker function, following an
original generalization of Whittaker models to incomplete flag
manifolds from \cite{GLO}; this integral formula has been announced
in \cite{GLO}. Our derivation involves a generalization of the
Gelfand-Zetlin realization to infinite-dimensional
$\CU(\frak{gl}_N)$-modules introduced in \cite{GKL}. Our main result
(Theorem 2.1) generalizes Theorem 1.1 of \cite{GLO} to arbitrary
Grassmannian $\Gr_{m,\,N}$. Moreover, our integral representation
verifies the conjectural integral formula from \cite{HV}, although
we construct another solution to the D-module with a different
asymptotic behavior.

The paper is organized as follows. In Section 1 we review on
parabolic Whittaker functions introduced in \cite{GLO}, and
formulate our main results: the Mellin-Barnes integral
representation for the specialized $\Gr_{m,\,N}$-Whittaker function
(Theorem 1.1), and its asymptotic behavior (Theorem 1.2). The second
part of the text contains a detailed proof of the main results. In
particular, we recall the generalized Gelfand-Zetlin realization of
the universal enveloping algebra $\CU(\gl_N)$ from \cite{GKL}, and
then we find out the Whittaker vectors (Proposition 2.1). In Section
3 we prove Theorems 1.1 and 1.2.

{\em Acknowledgments}: The author is thankful to A. Gerasimov and D.
Lebedev for very useful discussions.

\section{The $\Gr_{m,\,N}$-Whittaker function and its integral representation}

Let $\gl_N$ be the Lie algebra of $(N\times N)$ real matrices with
the Cartan subalgebra $\frak{h}\subset\gl_N$ of diagonal matrices,
and let $\frak{b}_{\pm}\subset\gl_N$ be a pair of opposed Borel
subalgebras containing $\frak{h}$. Then one has the triangular
decomposition $\gl_N=\frak{n}_-\oplus\frak{h}\oplus\frak{n}_+$,
where $\frak{n}_{\pm}\subset\frak{b}_{\pm}$ are the nilpotent
radicals given by strictly lower- and upper-triangular matrices. In
this way, the set of roots $R\subset\frak{h}^*$ decomposes into
$R_+\sqcup R_-$, where $R_+\subset R\subset\frak{h}^*$ is the set of
positive roots. Identifying $\frak{h}\simeq\IR^N$ with coordinates
$\ux=(x_1,\ldots,x_N)$ one may write
$R=\{\alpha\in\frak{h}^*|\,\alpha(\ux)=x_i-x_j,\,i\neq j\}$ and
$R_+=\{\alpha\in\frak{h}^*|\,\alpha(\ux)=x_i-x_j,\,i<j\}$. Clearly,
positive roots span the Borel subalgebra $\frak{b}_+$, and $R_-$
span $\frak{b}_-\subset\gl_N$. Let $\Delta\subset R_+$ be the set of
simple roots $\alpha_i(\ux)=x_i-x_{i+1}\in\frak{h}^*,\,1\leq i\leq
N-1$, and let $\{\omega_m,\,1\leq m\leq N\}$ be the non-reduced set
of fundamental weights given by $\omega_m(\ux)=x_1+\ldots+x_m$. The
Weyl group $\frak{S}_N$ is generated by simple reflections
$s_i=s_{\alpha_i}$, and acts in $\frak{h}^*$ by linear
transformations:
$$
 s_i(\beta)\,=\,\beta\,-\,(\alpha_i,\beta)\alpha_i\,,
\hspace{1.5cm}
 \beta\in\frak{h}^*.
$$
In particular one has $\frak{S}_N\cdot R_+=R_-$. Let
$I=\{1,2,\ldots,N-1\}$ be the set of vertices of Dynkin diagram,
then given a subset $J\subseteq I$, let $\ov{J}=I\setminus J$, and
let us consider the subgroup $\frak{W}_J\subset\frak{S}_N$ generated
by $\{s_j,\,j\in\ov{J}\}$. Then let $R_J\subseteq R_+$ be a subset
of positive roots defined by $\frak{W}_J\cdot R_J=-R_J$, and let
$\ov{R}_{J}=R_+\setminus R_J$. Then the corresponding parabolic
subalgbera is spanned by $R_-$ and $R_J$, and the corresponding
parabolic subgroup is denoted by $P_J$. In this paper we restrict
ourselves to the case $J=\{m\}\subset\{1,2,\ldots,N-1\}$ with
$\frak{W}_m=\frak{S}_m\times\frak{S}_{N-m}$, and $GL_N(\IC)/P_m$
being isomorphic to the Grassmannian $\Gr_{m,N}$. In this case we
have $I=I'\sqcup\,I''$ with $I'=\{1,\ldots,m\}$ and
$I''=\{m+1,\ldots,N\}$, then $\ov{R}_m$ is spanned by positive roots
$\alpha$ of the form $\alpha(\ux)=x_i-x_j;\,i\in I',\,j\in I''$.

Next, let us recall an original construction of
$\Gr_{m,N}$-Whittaker functions from \cite{GLO}. Let $B=B_-\subset
GL_N(\IC)$ be the Borel subgroup of lower-triangular matrices, and
let us pick a character $\chi_{\ul}:B_-\to\IC$ defined by
$\ul=(\la_1,\ldots,\la_N)\in\IC^N$. Then the associated Whittaker
function is defined as a certain matrix element of a principle
series representation $\CV_{\ul}=\Ind_B^{GL_N}\chi_{\ul}$.

Let us associate with $P_m$ a decomposition of the Borel subalgebra
$\frak{b}_+\subset\frak{gl}_N$
$$
 \frak{b}_+\,\,=\,\,\frak{h}^{(m)}
 \,\,\oplus\,\,
 \frak{n}_+^{(m)}\,,
$$
into the commutative subalgebra $\frak{h}^{(m)}\subset\frak{b}_+$
spanned by
 \be\label{Cartanml}
  H_1\,=\,E_{11}+\ldots+E_{mm}\,;
  \hspace{1.5cm}
  H_k\,=\,E_{1,k},\,\quad 2\leq k\leq m\,;\\
  H_{m+k}\,=\,E_{m+k,\,\ell+m},\,\quad
  1\leq k\leq N-m-1\,;\\
  H_{\ell+m}\,=\,E_{m+1,\,m+1}+\ldots+E_{\ell+m,\,\ell+m}\,,
 \ee
and the Lie subalgebra $\frak{n}^{(m)}_+\subset\frak{b}_+$ generated
by
 \be
  \frak{n}^{(m)}_+\,=\,\bigl\<E_{1,\,\ell+m};\,E_{1,\,m+1};\,
  E_{m,\,\ell+m};\hspace{5cm}\\
  \hspace{2cm}
  E_{kk}\,,2\leq k\leq N-1\,;\,
  E_{j,\,j+1},\,2\leq j\leq N-2\bigr\>\,.
 \ee
Note that $\dim\frak{h}^{(m)}=\rk\,\frak{gl}_N=N$ and
$\dim\frak{n}^{(m)}_+=N(N-1)/2$. Let $H^{(m)}$ and $N_+^{(m)}$ be
the  Lie groups corresponding to  the Lie algebras
$\mathfrak{h}^{(m)}$ and $\mathfrak{n}_+^{(m)}$. An open part
$GL_N^{\circ}$ (the big Bruhat cell) of $GL_N$ allows the following
analog of the Gauss decomposition:
 \be\label{modG}
  GL_N^{\circ}\,=\,N_-\,H^{(m)}\,N_+^{(m)}\,.
 \ee
Let $\CU=\CU(\gl_N)$ be the universal enveloping algebra of $\gl_N$.
The principal series representation $\CV_{\ul}$ admits a natural
structure of $\CU$-module, as well, as a module over the opposite
algebra $\CU^{\rm opp}$. Let us assume that the action of the Cartan
subalgebra $\frak{h}\subset\gl_N$ in $\CV_{\ul}$ is integrable to
the action of the Cartan torus $H\subset GL_N(\IR)$. Below we
introduce a pair of elements, $\<\psi_L|,\,|\psi_R\>\in\CV_{\ul}$,
generating a pair of dual submodules, $\CW_L=\<\psi_L|\CU^{\rm opp}$
and $\CW_R=\CU|\psi_R\>$, in $\CV_{\ul}$ (we adopt the bra- and ket-
vector notations to distinguish $\CU$- and $\CU^{\opp}$-modules
structures on $\CW_R$ and $\CW_L$ respectively).
\begin{de}\cite{GLO} The  $\Gr_{m,N}$-Whittaker vectors
$\<\psi_L|\in \CV'_{\underline{\lambda}}$ and
$|\psi_R\>\in\CV_{\underline{\lambda}}$ are defined by the following
conditions:
 \be\label{LeftWhittEqs}
  \<\psi_L|E_{n+1,\,n}\,=\,\h^{-1}\<\psi_L|\,,
  \hspace{2.5cm}
  1\leq n\leq N-1\,,
 \ee
 \be\label{RightWhittEqs}
  \left\{
  \begin{array}{lc}
  E_{kk}|\psi_R\>\,=\,0\,,&2\leq k\leq N-1\,;\\
  E_{k,\,k+1}|\psi_R\>\,=\,0,\,&2\leq k\leq N-2\,;\\
  E_{1,\,m+1}|\psi_R\>\,=E_{m,\,N}|\psi_R\>\,=\,0\,;&\\
  E_{1,\,N}|\psi_R\>\,=(-1)^{\epsilon(m,N)}\frac{1}{\h}|\psi_R\>&
  \end{array}
  \right.
 \ee
where  $\epsilon(m,N)$ is an integer number and $\hbar \in \IR$.
\end{de}

Note that the equations \eqref{LeftWhittEqs} define a
one-dimensional representation $\<\psi_L|$ of the Lie algebra
$\frak{n}_-$ of strictly lower-triangular matrices, and the
equations \eqref{RightWhittEqs} define a one-dimensional
representation of $\frak{n}^{(m)}_+$.
\begin{de}\cite{GLO} The $\Gr_{m,N}$-Whittaker function associated
with  the principal series representation
$\bigl(\pi_{\ul},\,\CV_{\ul}\bigr)$ is defined as the following
matrix element:
 \be\label{GrWhittaker}
  \Psi^{(m,N)}_{\underline{\lambda}}(\ux)\,
  =\,e^{-x_1\frac{m(N-m)}{2}}
  \bigl\<\psi_L\bigl|
  \pi_{\underline{\lambda}}\bigl(g(x_1,\ldots,x_N)\bigr)\bigr|
  \psi_R\bigr\>\,,
 \ee
where  the left and right  vectors  solve the equations
\eqref{LeftWhittEqs} and \eqref{RightWhittEqs} respectively. Here
$g(x)$ is a Cartan group valued function given by
 \be\label{GrGroupElement}
  g(\ux)\,=\,\exp\Big\{-\sum_{i=1}^N\,x_iH_i\Big\}\,,
 \ee
where $\ux=(x_1,\ldots,x_N)$ and the generators $H_i$,
$i=1,\ldots,N$ are defined by \eqref{Cartanml}.
\end{de}

In \cite{GLO} (Theorem 1.1) an  integral representation of the
$\Gr_{m,N}$-Whittaker function \eqref{GrWhittaker} was constructed
in the case of projective space $\IP^{N-1}$, corresponding to $m=1$.
We generalize this construction to generic Grassmannians
$\Gr_{m,N}$.
\begin{te} The specialized
$\Gr_{m,N}$-Whittaker function possesses the following integral
representation:
 \be\label{MBIntegralRep}
  \Psi^{(m,\,N)}_{\ul}(x,0,\ldots,0)\,=\,\int\limits_{\CC}\!
  d\ug\,\,e^{-\frac{x}{\h}\sum\limits_{i=1}^m\gamma_i}\,\,
  \frac{\prod\limits_{i=1}^m\prod\limits_{j=1}^N
  \Gamma_1\bigl(\gamma_i-\la_j|\h\bigr)}
  {\prod\limits_{i,k=1\atop k\neq i}^m
  \Gamma_1\bigl(\gamma_i-\gamma_k|\h\bigr)}\,,
 \ee
with $\underline{\gamma}=(\gamma_1,\ldots,\gamma_m)$ and
$\ul=(\la_1,\ldots,\la_N)\in\IR^N$. The integration contour is given
by $\CC=\bigl(\i\IR+\e\bigr)^m$, where $\e>\max\limits_{1\leq j\leq
N}\{\la_j\}$.
\end{te}
Here we use the following normalization of classical Gamma-function:
$$
 \Gamma_1\bigl(z\bigr|\h\bigr)\,\,
 =\,\,\h^{\frac{z}{\h}}\,\Gamma\Big(\frac{z}{\h}\Big)\,.
$$
We prove Theorem 1.1 in Section 3.

Clearly, the integral \eqref{MBIntegralRep} converges absolutely due
to the Stirling formula:
$$
 \Gamma\bigl(z+\la\bigr)\,
 =\,\sqrt{2\pi}\,z^{z+\la-\frac{1}{2}}\,e^{-z}\Big[1\,+\,O(z^{-1})\Big]\,,
\hspace{1.5cm}
 z\to\infty\,,
$$
when $|\arg(z)|<\pi$.

Integral representation \eqref{MBIntegralRep} coincides with the
expected one (5.1) in \cite{GLO}. Besides, a similar integral
formula was conjectured in \cite{HV} (see formulas (A.1) and (A.2)
in Appendix), but with different integration measure
$\widetilde{\mu}=\prod\limits_{i<j}(\ga_i-\ga_j)$. Our choice of
measure \eqref{GZMeasure} is provided by the generalized
Gelfand-Zetlin realization \cite{GKL}, and it is crucial in our
representation theory framework. Actually, the two solutions given
by integral formulas \eqref{MBIntegralRep} and the one from
\cite{HV}, have different asymptotic behavior, and below we derive
asymptotic of our solution.
\begin{te} When $x\to-\infty$, the specialized $\Gr_{m,N}$-Whittaker
function has the following asymptotic behavior:
 \be\label{Asymptotics}
  \Psi^{(m,\,N)}_{\underline{\gamma}_N}(x,0,\ldots,0)\,\,
  \sim\,\,m!\!\!
  \sum_{\s\in\,\frak{S}_N\!\bigl/\frak{W}_m}
  e^{-x\bigl(\s\cdot\,\omega_m(\ul)\bigr)}
  \bigl(\s\cdot c_m\bigr)(\ul)\,.
 \ee
with $\omega_m(\ul)=\la_1+\ldots+\la_m$, and
 \be
  c_m(\ul)\,=\,\prod_{\alpha\in\ov{R}_m}
  \Gamma_1\bigl(\alpha(\ul)\bigr|\h\bigr)\,,
 \hspace{1.5cm}
  \bigl(\s\cdot c_m\bigr)(\ul)\,
  =\,\prod_{\alpha\in\ov{R}_m}
  \Gamma_1\bigl(\s\cdot\alpha(\ul)\bigr|\h\bigr)\,.
 \ee
\end{te}
A proof of Theorem 1.2 is given in Section 3.1.

\section{Construction of $\Gr_{m,N}$-Whittaker vectors}

In this Section we construct explicit solutions to
\eqref{LeftWhittEqs} and \eqref{RightWhittEqs} using the generalized
Gelfand-Zetlin realization of principal series
$U(\frak{gl}_N)$-modules from \cite{GKL}. Namely, let
$\ug_1,\ldots,\ug_N$ be a triangular array consisting of $N(N-1)/2$
variables
$\underline{\gamma}_n=(\gamma_{n1},\ldots,\gamma_{nn})\in\IC^n,n=1,\ldots,N$.
The following operators define a representation $\pi$ of $\CU$ in
the space $\CM_N$ of meromorphic functions in $N(N-1)/2$ variables
$(\ug_1,\ldots,\ug_{N-1})$:
 \be\label{GZRep}
  E_{kk}\,=\,\frac{1}{\h}\Big(\sum_{j=1}^n\gamma_{n,j}\,
  -\,\sum_{i=1}^{n-1}\gamma_{n-1,\,i}\Big)\,,
  \hspace{2.5cm}
  1\leq k\leq N\,;\\
  E_{n,\,n+1}\,
  =\,-\frac{1}{\h}\sum_{i=1}^n\frac{\prod\limits_{j=1}^{n+1}
  (\gamma_{n,i}-\gamma_{n+1,\,j}-\frac{\h}{2})}
  {\prod\limits_{s\neq i}(\gamma_{n,i}-\gamma_{n,s})}\,
  e^{-\h\pr_{n,i}}\,,
  \hspace{1.5cm}
  1\leq n\leq N-1\,;\\
  E_{n+1,\,n}\,
  =\,\frac{1}{\h}\sum_{i=1}^n\frac{\prod\limits_{j=1}^{n-1}
  (\gamma_{n,i}-\gamma_{n-1,\,j}+\frac{\h}{2})}
  {\prod\limits_{s\neq i}(\gamma_{n,i}-\gamma_{n,s})}\,
  e^{\h\pr_{n,i}}\,,
  \hspace{1.5cm}
  1\leq n\leq N-1\,,
 \ee
where $E_{ij}=\pi(e_{ij}),\,1\leq i,j=1\leq N$ for the standard
elementary matrix units $e_{ij}\in\gl_N$. This realization of
universal enveloping algebra $\CU$ is referred to as generalized
Gelfand-Zetlin realization.
\begin{rem} Evidently the Weyl group $\frak{S}_N$ acts on $E_{ij}$
in \eqref{GZRep} by permutations of indices $(i,j)$. This provides
$N!$ different realizations of $\CU$, and we use certain
$\frak{S}_N$-twisted generalized Gelfand-Zetlin realizations of
$\CU$ in the next Section for derivation of Whittaker vectors
\eqref{LeftWhittEqs}, \eqref{RightWhittEqs}.
\end{rem}

The the universal enveloping algebra $\CU$ acts in
$\CW_R\subseteq\CM_N$ by differential operators\eqref{GZRep}, and
the opposite algebra $\CU^{\opp}$ acts in $\CW_L\subseteq\CM_N$ via
the adjoint operators:
 \be
  E^{\dag}_{ij}\,=\,\mu(\ga)^{-1}\,E_{ij}\,\mu(\ga)\,,
 \hspace{1.5cm}
  1\leq i,j\leq N\,,
 \ee
with
 \be\label{GZMeasure}
  \mu(\ga)\,=\,\prod_{n=2}^{N-1}\mu_n(\ug_n)\,
  =\,\prod_{n=2}^{N-1}\prod_{i,j=1\atop i\neq j}^n\,
  \frac{1}{\Gamma\bigl(\frac{\ga_{ni}-\ga_{nj}}{\h}\bigr)}\,.
 \ee
It was shown in \cite{GKL} that there exist a non-degenerate pairing
between the modules $\CW_L$ and $\CW_R$ with measure
\eqref{GZMeasure}:
 \be\label{Pairing}
  \<\phi_1,\,\phi_2\>\,=\,\int\limits_{\IR^{N(N-1)/2}}
  \phi_1(\ga)\phi_2(\ga)\,\,
  \mu(\ga)\prod_{k,n=1\atop k\leq n}^Nd\ga_{nk}\,,
 \ee
where $\phi_1\in\CW_L$ and $\phi_2\in\CW_R$.
\begin{prop} For $1<m<N$ the Whittaker vectors have the following
expressions.
\begin{enumerate}
\item  A solution to \eqref{LeftWhittEqs} is given by
 \be\label{LeftWhitVec}
  \psi_L^{(m)}\,=\,e^{\i\pi\gamma_{11}}\prod_{i=1}^{m-1}\prod_{j=1}^m
  \frac{1}{\Gamma_1(\gamma_{m-1,\,i}-\gamma_{m,j}+\frac{\h}{2}|\h)}\,.
 \ee
\item A solution to \eqref{RightWhittEqs} is given by
 \be\label{RightWhitVec}
  \psi_R^{(m)}\,=\,\delta(\gamma_{11})\prod_{i=1}^m
  \prod_{j=1}^N\Gamma_1\Big(\gamma_{N-1,\,i}-\gamma_{Nj}+\frac{\h}{2}
  \Big|\h\Big)\cdot
  \prod_{a=1}^{m-1}\prod_{b=1}^m
  \Gamma_1\Big(\gamma_{m-1,\,a}-\gamma_{mb}+\frac{\h}{2}\Big|\h\Big)\\
  \times\prod_{n=2\atop n\neq m}^{N-1}
  \Big[\delta\Big(\sum_{j=1}^n\gamma_{nj}
  -\sum_{i=1}^{n-1}\gamma_{n-1,\,i}\Big)
  \prod_{k=1}^{n-1}
  \delta\Big(\gamma_{n-1,\,k}-\gamma_{nk}+\frac{\h}{2}\Big)
  \prod_{i,j=1\atop i\neq j}^n\Gamma_1\bigl(\gamma_{ni}-\gamma_{nj}
  \bigl|\h\bigr)\Big]
  \,,
 \ee
\end{enumerate}
\end{prop}
Further, substituting the Whittaker vectors \eqref{LeftWhitVec} and
\eqref{RightWhitVec} into the pairing \eqref{Pairing} we arrive to
the following integral representation.

\subsection{Proof of Proposition}

Due to the action of the Weyl group $\frak{S}_N$, actually one has
$N!$ realizations of $U(\frak{gl}_N)$ defined by
 \be\label{TwistedGZRep}
  E_{ij}^w\,:=\,wE_{ij}w^{-1}\,,
 \hspace{2.5cm}
  w\in\frak{S}_N\,,
 \hspace{1cm}
  i,j=1,\ldots,N\,.
 \ee
Let us call these realizations of $U(\frak{gl}_N)$ the $w$-twisted
Gelfand-Zetlin realizations.

Given the simple reflections $s_i\in\frak{S}_N,\,i\in I$, let us
introduce the Coxeter elements $c_n=s_n\cdot\ldots\cdot s_1,\,n\in
I$. In particular, one has $c_1=s_1$, and for the longest element
$w_0\in\frak{S}_N$ the following decomposition holds:
$$
 w_0\,=\,c_1c_2\cdot\ldots\cdot c_{N-1}\,.
$$

\emph{Proof of Proposition 2.1.} Given $m>1$ let us solve the
defining relations \eqref{LeftWhittEqs}, \eqref{RightWhittEqs},
using the $c_{m-1}$-twisted Gelfand-Zetlin realization
\eqref{TwistedGZRep}.

\npaa Let us start from solving the equations \eqref{LeftWhittEqs}
for the left $\Gr_{m,N}$-Whittaker vector. Namely, one has to check
that \eqref{LeftWhitVec} satisfies \eqref{LeftWhittEqs}:
 \be
  \bigl(E_{k+1,\,k}^{c_{m-1}}\bigr)^{\dag}\psi_L^{(m)}\,
  =\,\h^{-1}\psi_L^{(m)}\,,
 \hspace{2.5cm}
  1\leq k<N\,.
 \ee
Actually, one has to check only two relations:
 \be\label{TwistedLeftWhittEqs}
  -\bigl(E_{21}^{c_{m-1}}\bigr)^{\dag}\psi_L^{(m)}\,
  =\,\bigl(E_{m+1,\,m}^{c_{m-1}}\bigr)^{\dag}\psi_L^{(m)}\,
  =\,\h^{-1}\psi_L^{(m)}\,,
 \ee
and the other relations are evidently true, since the difference
operators $E_{k+1,\,k}^{c_{m-1}},\,k\neq 1,m$ act trivially on
$\psi_L^{(m)}=\psi_L^{(m)}\bigl(\underline{\gamma}_{m-1},\,\underline{\gamma}_m\bigr)$.
The relations \eqref{TwistedLeftWhittEqs} can be verified using the
following combinatorial formulas.

\begin{lem}
Given a set of variables
$\underline{\gamma}=(\gamma_1,\ldots,\gamma_n)$ the following
identities hold:
\begin{enumerate}
\item
 \be\label{Combin1}
  \sum_{i=1}^n\,\gamma_i^m\prod_{i\neq k}\frac{1}{\gamma_i-\gamma_k}\,
  =\,\delta_{m,\,n-1}\,,
 \hspace{1.5cm}
  m<n\,;
 \ee
More generally, one has
$$
 \sum_{i=1}^n\,\gamma_i^m\prod_{i\neq k}\frac{1}{\gamma_i-\gamma_k}\,
 =\,\sum_{k_1+\cdots+k_n=n+1-m}
 \gamma_1^{k_1}\cdot\ldots\cdot\gamma_n^{k_n}\,.
$$
\item
 \be\label{Combin2}
  \sum_{i=1}^n\,\prod_{i\neq k}\frac{c-\gamma_k}{\gamma_i-\gamma_k}\,
  =\,1\,,
 \ee
for any constant $c$.
\end{enumerate}
\end{lem}
\emph{Proof}. We have
$$
 \sum_{i=1}^n\gamma_i^m\prod_{i\neq k}\frac{1}{\gamma_i-\gamma_k}\,=
 \oint\limits_{\gamma_i}d\lambda\frac{\lambda^m}{A_n(\lambda)}\,,
 \hspace{1.5cm}
 A_n(\lambda)=\prod_{i=1}^n(\lambda-\gamma_i)\,.
$$
Taking the residue at infinity we get
$$
 \oint_{\infty}d\lambda\,\frac{\lambda^{m-n}}
 {1+\sum\limits_{k=1}^n(-1)^k\s_k(\underline{\gamma})\lambda^{-k}}\,
 =\,\chi_{n+1-m}(\underline{\gamma})\,\Theta(m+1-n)\,,
$$
where
$$
 \s_k(\underline{\gamma})\,
 =\,\sum_{i_1<\ldots<i_k}\gamma_{i_1}\cdot\ldots\cdot\gamma_{i_k}\,,
\hspace{1.5cm}
 \chi_k(\underline{\gamma})\,
 =\,\sum_{i_1+\ldots+i_n=k}\gamma_1^{k_1}\cdot\ldots\cdot\gamma_n^{k_n}
$$
are the characters of finite-dimensional representations
$\bigwedge^k\IC^n$ and ${\rm Sym}^k\IC^n$ respectively.

The other identity can be proved similarly. $\Box$

Then for the second relation in \eqref{TwistedLeftWhittEqs} one
readily derives the following:
 \be\label{InductiveStep}
 \hspace{-5cm}
  \bigl(E_{m+1,\,m}^{c_{m-1}}\bigr)^{\dag}\psi_L^{(m)}\,
  =\,E_{m+1,\,m-1}^{\dag}\psi_L^{(m)}\,
  =\,\frac{1}{\h}\sum_{i_1=1}^m\prod_{k_1=1\atop k_1\neq i_1}^m
  \frac{1}{\gamma_{m,\,i_1}-\gamma_{m,\,k_1}}\\
  \times\sum_{i_2=1}^{m-1}\prod_{k_2=1\atop k_2\neq i_2}^{m-1}
  \frac{(\gamma_{m,\,i_1}-\gamma_{m-1,\,k_2}-\frac{\h}{2})}
  {\gamma_{m-1,\,i_2}-\gamma_{m-1,\,k_2}}\,
  \prod_{j=1}^{m-2}\bigl(\gamma_{m-1,\,i_2}
  -\gamma_{m-2,\,j}-\frac{\h}{2}\bigr)\,
  e^{-\h(\pr_{m,\,i_1}+\pr_{m-1,\,i_2})}\cdot\psi_L^{(m)}\\
 \hspace{-1.5cm}
  =\,\frac{1}{\h}\sum_{i_1=1}^m\prod_{k_1=1\atop k_1\neq i_1}^m
  \frac{1}{\gamma_{m,\,i_1}-\gamma_{m,\,k_1}}
  \sum_{i_2=1}^{m-1}
  \frac{\prod\limits_{j_1=1\atop j_1\neq i_1}^m
  \bigl(\gamma_{m-1,\,i_1}-\gamma_{m,\,j_1}-\frac{\h}{2}\bigr)}
  {\prod\limits_{k_2\neq i_2}
  \bigl(\gamma_{m-1,\,i_2}-\gamma_{m-1,\,k_2}\bigr)}\,
  \prod_{j_2=1}^{m-2}
  \bigl(\gamma_{m-1,\,i_2}-\gamma_{m-2,\,j}-\frac{\h}{2}\bigr)
  \psi_L^{(m)}\,=\ldots
 \ee
Using \eqref{Combin1} we have
$$
 \ldots=\,\frac{1}{\h}\sum_{i_1=1}^m\prod_{k_1=1\atop k_1\neq i_1}^m
 \frac{\gamma_{m,\,k_1}}{\gamma_{m,\,i_1}-\gamma_{m,\,k_1}}
 \sum_{i_2=1}^{m-1}
 \frac{-\gamma_{m-1,\,i_2}^{m-2}}
 {\prod\limits_{k_2\neq i_2}
 \bigl(\gamma_{m-1,\,i_2}-\gamma_{m-1,\,k_2}\bigr)}\psi_L^{(m)}\,
 =\ldots
$$
and finally we apply \eqref{Combin2} and obtain
$$
 =\,\frac{1}{\h}\sum_{i_1=1}^m\prod_{k_1=1\atop k_1\neq i_1}^m
 \frac{-\gamma_{m,\,k_1}}{\gamma_{m,\,i_1}-\gamma_{m,\,k_1}}\,
 =\,\frac{1}{\h}\psi_L^{(m)}\,.
$$

The first relation in \eqref{TwistedLeftWhittEqs} reads as follows:
 \be\label{E1mIdentity}
 \hspace{-2.5cm}
  \bigl(E_{21}^{c_{m-1}}\bigr)^{\dag}\psi_L^{(m)}\,
  =\,E_{1m}^{\dag}\psi_L^{(m)}\,
  =\,-\frac{1}{\h}\sum_{i_1=1}^{m-1}
  \frac{\prod\limits_{j_1=1}^{m+1}\bigl(\gamma_{m-1,\,i_1}
  -\gamma_{m,\,j_1}+\frac{\h}{2}\bigr)}
  {\prod\limits_{k_1\neq i_1}\bigl(\gamma_{m-1,\,i_1}-\gamma_{m-1,\,k_1}\bigr)}
  \sum_{i_2=1}^{m-2}
  \frac{\prod\limits_{j_2=1\atop j_2\neq i_1}^m\bigl(\gamma_{m-2,\,i_2}
  -\gamma_{m-1,\,j_2}+\frac{\h}{2}\bigr)}
  {\prod\limits_{k_2\neq i_2}\bigl(\gamma_{m-2,\,i_2}-\gamma_{m-2,\,k_2}\bigr)}
  \times\ldots\\
 \hspace{-1.5cm}
  \times
  \sum_{i_{m-2}=1}^2
  \frac{\prod\limits_{j_{m-2}=1\atop j_{m-2}\neq i_{m-1}}^3\!\!\!
  \bigl(\gamma_{2,\,i_{m-2}}-\gamma_{3,\,j_{m-2}}+\frac{\h}{2}\bigr)}
  {\prod\limits_{k_{m-2}\neq i_{m-2}}\!\!\!
  \bigl(\gamma_{2,\,i_{m-2}}-\gamma_{2,\,k_{m-2}}\bigr)}
 \hspace{-5mm}
  \prod_{j_{m-1}=1\atop j_{m-1}\neq i_{m-2}}^2
 \hspace{-3mm}
  \Big(\gamma_{11}-\gamma_{2,\,j_{m-1}}+\frac{\h}{2}\bigr)\,
  e^{^{\h\bigl(\pr_{11}+\sum\limits_{n=2}^{m-1}\pr_{n,\,m-n}\bigr)}}\!\!
  \cdot\psi_L^{(m)}\,\,
  =\,\,-\frac{1}{\h}\psi_L^{(m)}
 \ee
Thus we have to check the following identity:
 \be
 \hspace{-1cm}
  \frac{1}{\h}\sum_{i_1=1}^{m-1}\prod\limits_{k_1\neq i_1}
  \frac{1}{\gamma_{m-1,\,i_1}-\gamma_{m-1,\,k_1}}
  \sum_{i_2=1}^{m-2}
  \frac{\prod\limits_{j_2=1\atop j_2\neq i_1}^m\bigl(\gamma_{m-2,\,i_2}
  -\gamma_{m-1,\,j_2}+\frac{\h}{2}\bigr)}
  {\prod\limits_{k_2\neq i_2}\bigl(\gamma_{m-2,\,i_2}-\gamma_{m-2,\,k_2}\bigr)}
  \times\ldots\\
  \times
  \sum_{i_{m-2}=1}^2
  \frac{\prod\limits_{j_{m-2}=1\atop j_{m-2}\neq i_{m-1}}^3\bigl(\gamma_{2,\,i_{m-2}}
  -\gamma_{3,\,j_{m-2}}+\frac{\h}{2}\bigr)}
  {\prod\limits_{k_{m-2}\neq i_{m-2}}
  \bigl(\gamma_{2,\,i_{m-2}}-\gamma_{2,\,k_{m-2}}\bigr)}
 \hspace{-5mm}
  \prod_{j_{m-1}=1\atop j_{m-1}\neq i_{m-2}}^2
 \hspace{-3mm}
  \Big(\gamma_{11}-\gamma_{2,\,j_{m-1}}+\frac{\h}{2}\bigr)\,
  =\,\frac{1}{\h}\,.
 \ee
This identity can be verified by induction over $m$. Indeed, for
$m=2$ \eqref{E1mIdentity} reads
 \be
 \hspace{-1.5cm}
  \bigl(E_{21}^{c_2}\bigr)^{\dag}\psi_L^{(2)}\,
  =\,E_{12}^{\dag}\psi_L^{(2)}\\
  =\,-\frac{1}{\h}\Big(\gamma_{11}-\gamma_{21}+\frac{\h}{2}\Big)
  \Big(\gamma_{11}-\gamma_{22}+\frac{\h}{2}\Big)e^{\h\pr_{11}}\cdot
  e^{\i\pi\gamma_{11}}
  \prod_{i=1}^2\frac{1}{\Gamma_1\bigl(\gamma_{11}-\gamma_{2i}+\frac{1}{2}
  \bigl|\h\bigr)}\\
  =\,\frac{1}{\h}\psi_L^{(2)}\,.
 \ee
The inductive step directly follows the reasoning from
\eqref{InductiveStep}, using the combinatorial identities
\eqref{Combin1} and \eqref{Combin2}.

\npaa Let us check that the expression \eqref{RightWhitVec}
satisfies the relations the $c_{m-1}$-twisted relations
\eqref{RightWhittEqs}:
 \be\label{TwistedRightWhittEqs}
  E_{nn}^{c_{m-1}}\psi_R^{(m)}\,=\,0,
 \hspace{5mm}
  n=2,\ldots,N-1\,;
 \hspace{1cm}
  E_{k,\,k+1}^{c_{m-1}}\psi_R^{(m)}\,=\,0,
 \hspace{5mm}
  k=2,\ldots,N-2\,;\\
  E_{1,\,m+1}^{c_{m-1}}\psi_R^{(m)}\,
  =\,E_{m,\,N}^{c_{m-1}}\psi_R^{(m)}\,=\,0\,;
 \hspace{1.5cm}
  E_{1,\,N}^{c_{m-1}}\psi_R^{(m)}\,
  =\,(-1)^{\epsilon(m,N)}\h^{-1}\psi_R^{(m)}\,.
 \ee
The relations corresponding to Cartan generators
$E_{nn}^{c_{m-1}}\psi_R^{(m)}=0\,,n=2,\ldots,N-1$ hold due to the
delta-factors
$$
 \delta(\gamma_{11})\prod_{n=2\atop n\neq m}^{N-1}
 \delta\Big(\sum_{j=1}^n\gamma_{nj}
 -\sum_{i=1}^{n-1}\gamma_{n-1,\,i}\Big)
$$
in \eqref{RightWhitVec}, and since
$$
 E_{kk}^{c_m}\,=\,E_{k-1,\,k-1},\,
 \hspace{5mm}
 k=2,\ldots,m\,;
 \hspace{1.5cm}
 E_{kk}^{c_m}\,=\,E_{kk},\,
 \hspace{5mm}
 k=m+1,\ldots,N-1\,.
$$
Similarly, the relations
$E_{k,\,k+1}^{c_{m-1}}\psi_R^{(m)}=0\,,k=2,\ldots,N-2$ hold due to
the delta-factors
 \be\label{DeltaFactor}
  \prod_{n=2\atop n\neq m}^{N-1}\prod_{k=1}^{n-1}
  \delta\Big(\gamma_{n-1,\,k}-\gamma_{nk}+\frac{\h}{2}\Big)
 \ee
in \eqref{RightWhitVec}, and since
$$
 E_{k,\,k+1}^{c_{m-1}}\,=\,E_{k-1,\,k},\,
 \hspace{5mm}
 k=2,\ldots,m\,;
 \hspace{1.5cm}
 E_{k,\,k+1}^{c_{m-1}}\,=\,E_{k,\,k+1},\,
 \hspace{5mm}
 k=m+1,\ldots,N-2\,.
$$
Due to the same delta-factor \eqref{DeltaFactor} one has
$E_{1,\,m+1}^{c_{m-1}}\psi_R^{(m)}\,=\,E_{m,\,m+1}\psi_R^{(m)}=0$.

Thus we have to check the remaining two relations:
 \be\label{TwistedRightWhittEqs2}
 \hspace{-1cm}
  E_{m,\,N}^{c_{m-1}}\psi_R^{(m)}\,=\,E_{m-1,\,N}\psi_R^{(m)}\,=\,0\,,
 \hspace{1cm}
  E_{1,\,N}^{c_{m-1}}\psi_R^{(m)}\,
  =\,E_{m,\,N}\psi_R^{(m)}\,
  =\,\frac{(-1)^m}{\h}\,\psi_R^{(m)}\,.
 \ee
\begin{lem}
For $n=1,\ldots,N-1$ the following holds:
 \be\label{EnNGenerator}
 \hspace{-5mm}
  E_{n,\,N}\,=\,-\frac{1}{\h}\sum_{i_1=1}^{N-1}
  \frac{\prod\limits_{j_1=1}^N\bigl(\gamma_{N-1,\,i_1}
  -\gamma_{N,\,j_1}-\frac{\h}{2}\bigr)}
  {\prod\limits_{k_1\neq i_1}\bigl(\gamma_{N-1,\,i_1}
  -\gamma_{N-1,\,k_1}\bigr)}\,\,\,
  \sum_{i_2=1}^{N-2}
  \frac{\prod\limits_{j_2=1\atop j_2\neq i_1}^{N-1}\bigl(\gamma_{N-2,\,i_2}
  -\gamma_{N-1,\,j_2}-\frac{\h}{2}\bigr)}
  {\prod\limits_{k_2\neq i_2}\bigl(\gamma_{N-2,\,i_2}-\gamma_{N-2,\,k_2}\bigr)}
  \times\ldots\\
 \hspace{-1.5cm}
  \times
  \sum_{i_{N-n}=1}^n\!\!\!
  \frac{\prod\limits_{j_{N-n}=1\atop j_{N-n}\neq i_{N+1-n}}^{N+1-n}\!\!\!
  \bigl(\gamma_{n,\,i_{N-n}}-\gamma_{n+1,\,j_{N-n}}-\frac{\h}{2}\bigr)}
  {\prod\limits_{k_{N-n}\neq i_{N-n}}\!\!\!
  \bigl(\gamma_{n,\,i_{N-n}}-\gamma_{n,\,k_{N-n}}\bigr)}\,\,
  e^{^{-\h\sum\limits_{a=n}^{N-1}\pr_{a,\,i_{N-a}}}}\,\,.
 \ee
\end{lem}
\emph{Proof.} Direct calculation using \eqref{GZRep}. $\Box$

At first let us note that due to the delta-factors
$$
 \prod_{n=m+1}^{N-1}\prod_{k=1}^n\delta\Big(\gamma_{n-1,i}
 -\gamma_{n,i}+\frac{\h}{2}\Big)
$$
in \eqref{RightWhitVec} one gets only $m$ non-vanishing terms in
\eqref{EnNGenerator} (for $n=m$):
 \be\label{EmNPsiRrelation1}
  E_{1,\,N}^{c_{m-1}}\psi_R^{(m)}\,
  =\,E_{m,\,N}\psi_R^{(m)}\\
  =\,-\frac{1}{\h}\Big\{\sum_{i_1=1}^{N-1}
  \frac{\prod\limits_{j_1=1}^N\bigl(\gamma_{N-1,\,i_1}
  -\gamma_{N,\,j_1}-\frac{\h}{2}\bigr)}
  {\prod\limits_{k_1\neq i_1}\bigl(\gamma_{N-1,\,i_1}
  -\gamma_{N-1,\,k_1}\bigr)}\,\,\,
  \sum_{i_2=1}^{N-2}
  \frac{\prod\limits_{j_2=1\atop j_2\neq i_1}^{N-1}\bigl(\gamma_{N-2,\,i_2}
  -\gamma_{N-1,\,j_2}-\frac{\h}{2}\bigr)}
  {\prod\limits_{k_2\neq i_2}\bigl(\gamma_{N-2,\,i_2}-\gamma_{N-2,\,k_2}\bigr)}
  \times\ldots\\
 \hspace{-1.5cm}
  \times
  \sum_{i_{N-m}=1}^m\,
  \frac{\prod\limits_{j_{N-m}=1\atop j_{N-m}\neq i_{N-m-1}}^m
  \bigl(\gamma_{m,\,i_{N-m}}-\gamma_{m+1,\,j_{N-m}}-\frac{\h}{2}\bigr)}
  {\prod\limits_{k_{N-m}\neq i_{N-m}}\!\!\!\!\!
  \bigl(\gamma_{m,\,i_{N-m}}-\gamma_{m,\,k_{N-m}}\bigr)}\,\,
  e^{^{-\h\sum\limits_{a=m}^{N-1}\pr_{a,\,i_{N-a}}}}\Big\}\cdot\psi_R^{(m)}\\
  =\,-\frac{1}{\h}\Big\{\sum_{i_1=1}^m
  \frac{\prod\limits_{j_1=1}^N\bigl(\gamma_{N-1,\,i_1}
  -\gamma_{N,\,j_1}-\frac{\h}{2}\bigr)}
  {\prod\limits_{k_1=1\atop k_1\neq i_1}^{N-1}\bigl(\gamma_{N-1,\,i_1}
  -\gamma_{N-1,\,k_1}\bigr)}\,
  \frac{\prod\limits_{j_2=1\atop j_2\neq i_1}^{N-1}\bigl(\gamma_{N-2,\,i_1}
  -\gamma_{N-1,\,j_2}-\frac{\h}{2}\bigr)}
  {\prod\limits_{k_2=1\atop k_2\neq i_1}^{N-1}\bigl(\gamma_{N-2,\,i_1}
  -\gamma_{N-2,\,k_2}\bigr)}\,\times\ldots\\
  \ldots\times
  \frac{\prod\limits_{j_{N-m}=1\atop j_{N-m}\neq i_1}^{m+1}
  \bigl(\gamma_{m,\,i_1}-\gamma_{m+1,\,j_{N-m}}-\frac{\h}{2}\bigr)}
  {\prod\limits_{k_{N-m}=1\atop k_{N-m}\neq i_1}^m\!\!\!\!\!
  \bigl(\gamma_{m,\,i_1}-\gamma_{m,\,k_{N-m}}\bigr)}\,\,
  e^{-\h\sum\limits_{a=m}^{N-1}\pr_{a,\,i_1}}\Big\}\cdot\psi_R^{(m)}\,=\ldots
 \ee
Secondly, taking into account that
 \be\label{ShiftGammaRelation1}
  e^{-\h\pr_{N-1,\,i}}\cdot
  \prod_{a=1}^m\prod_{b=1}^N\Gamma_1\Big(\gamma_{N-1,\,a}
  -\gamma_{N,\,b}+\frac{\h}{2}\Big|\,\h\Big)\\
  =\,\prod_{j=1}^N\frac{1}{\gamma_{N-1,\,i}-\gamma_{N,\,j}-\frac{\h}{2}}\,
  \prod_{a=1}^m\prod_{b=1}^N\Gamma_1\Big(\gamma_{N-1,\,a}
  -\gamma_{N,\,b}+\frac{\h}{2}\Big|\,\h\Big)\cdot
  e^{-\h\pr_{N-1,\,i}}\,;
 \ee
and due to the factors
 $\prod\limits_{n=m+1}^{N-1}\,\prod\limits_{i,j=1\atop i\neq j}^n
 \Gamma_1\bigl(\gamma_{ni}-\gamma_{nj}\bigl|\h\bigr)$ in \eqref{RightWhitVec} one has
 \be\label{EmNPsiRrelation2}
 \hspace{-1cm}
  \ldots=\,-\frac{1}{\h}\Big\{\sum_{i_1=1}^m
  \prod\limits_{k_1=1\atop k_1\neq i_1}^{N-1}
  \frac{1}{\bigl(\gamma_{N-1,\,i_1}
  -\gamma_{N-1,\,k_1}-\h\bigr)}\,
  \frac{\prod\limits_{j_2=1\atop j_2\neq i_1}^{N-1}\bigl(\gamma_{N-2,\,i_1}
  -\gamma_{N-1,\,j_2}-\frac{\h}{2}\bigr)}
  {\prod\limits_{k_2=1\atop k_2\neq i_1}^{N-1}\bigl(\gamma_{N-2,\,i_q}
  -\gamma_{N-2,\,k_2}-\h\bigr)}\,\times\ldots\\
 \hspace{-1.5cm}
  \ldots\times
  \frac{\prod\limits_{j_{N-m-1}=1\atop j_{N-m-1}\neq i_1}^{m+2}
  \bigl(\gamma_{m+1,\,i_1}-\gamma_{m+2,\,j_{N-m-1}}-\frac{\h}{2}\bigr)}
  {\prod\limits_{k_{N-m-1}=1\atop k_{N-m-1}\neq i_1}^{m+1}\!\!\!\!\!
  \bigl(\gamma_{m+1,\,i_1}-\gamma_{m+1,\,k_{N-m-1}}-\h\bigr)}\,\,
  \frac{\prod\limits_{j_{N-m}=1\atop j_{N-m}\neq i_1}^{m+1}
  \bigl(\gamma_{m,\,i_1}-\gamma_{m+1,\,j_{N-m}}-\frac{\h}{2}\bigr)}
  {\prod\limits_{k_{N-m}=1\atop k_{N-m}\neq i_1}^m\!\!\!\!\!
  \bigl(\gamma_{m,\,i_1}-\gamma_{m,\,k_{N-m}}\bigr)}\\
  \times\,\,
  e^{-\h\sum\limits_{a=m}^{N-1}\pr_{a,\,i_1}}\Big\}\cdot\psi_R^{(m)}\,
  =\ldots
 \ee
Next, since
 \be\label{DeltaMeasureRelation}
  \prod_{j_{a+1}=1\atop j_{a+1}\neq i_{a+1}}^{N-a}
 \hspace{-4mm}
  \Big(\gamma_{N-a-1,\,i_{a+1}}
  -\gamma_{N-a,\,j_{a+1}}-\frac{\h}{2}\Big)\equiv\\
  \equiv\,
  \prod_{k_a\neq i_a}\!\!
  \bigl(\gamma_{N-a,\,i_a}-\gamma_{N-a,\,k_a}-\h\bigr)\,\,\,
  \mod
  \prod_{n=m+1}^{N-1}
  \prod_{k=1}^{n-1}
  \delta\Big(\gamma_{n-1,\,k}-\gamma_{nk}-\frac{\h}{2}\Big)
 \ee
and since the factor $\prod\limits_{a=1}^{m-1}\prod\limits_{b=1}^m
  \Gamma_1\Big(\gamma_{m-1,\,a}-\gamma_{mb}+\frac{\h}{2}\Big|\h\Big)$
in \eqref{RightWhitVec} produces
 \be\label{ShiftGammaRelation2}
  e^{-\h\pr_{m,\,i}}\cdot
  \prod_{a=1}^{m-1}\prod_{b=1}^m\Gamma_1\Big(\gamma_{m-1,\,a}
  -\gamma_{n,\,b}+\frac{\h}{2}\Big|\,\h\Big)\\
  =\,\prod\limits_{r=1}^{m-1}\!\bigl(\gamma_{m-1,\,r}
  -\gamma_{m,\,i}-\frac{\h}{2}\bigr)\,\,\,
  \prod_{a=1}^{m-1}\prod_{b=1}^m\Gamma_1\Big(\gamma_{m-1,\,a}
  -\gamma_{m,\,b}+\frac{\h}{2}\Big|\,\h\Big)\cdot
  e^{-\h\pr_{m,\,i}}\,,
 \ee
one arrives to the following:
 \be\label{EmNPsiRrelation3}
  \ldots=\,-\frac{1}{\h}\sum_{i_1=1}^m\,
  \frac{\prod\limits_{r=1}^{m-1}
  \bigl(\gamma_{m-1,\,r}-\gamma_{m,\,i_1}-\frac{\h}{2}\bigr)}
  {\prod\limits_{k_{N-m}=1\atop k_{N-m}\neq i_1}^m\!\!
  \bigl(\gamma_{m,\,i_1}-\gamma_{m,\,k_{N-m}}\bigr)}\,\,\psi_R^{(m)}\,
  =\,\frac{(-1)^m}{\h}\psi_R^{(m)}\,,
 \ee
where the last equality follows from \eqref{Combin1}.

At last we have to verify the remaining first relation in
\eqref{TwistedRightWhittEqs2}. Following the same reasoning as in
\eqref{EmNPsiRrelation1}-\eqref{EmNPsiRrelation2} above, we obtain
the following:
 \be
  E_{m,\,N}^{c_{m-1}}\psi_R^{(m)}\,=\,E_{m-1,\,N}\psi_R^{(m)}\\
  =\,-\frac{1}{\h}\Big\{\sum_{i_1=1}^m
  \prod\limits_{k_1=1\atop k_1\neq i_1}^{N-1}
  \frac{1}{\bigl(\gamma_{N-1,\,i_1}
  -\gamma_{N-1,\,k_1}-\h\bigr)}\,
  \frac{\prod\limits_{j_2=1\atop j_2\neq i_1}^{N-1}\bigl(\gamma_{N-2,\,i_1}
  -\gamma_{N-1,\,j_2}-\frac{\h}{2}\bigr)}
  {\prod\limits_{k_2=1\atop k_2\neq i_1}^{N-1}\bigl(\gamma_{N-2,\,i_1}
  -\gamma_{N-2,\,k_2}-\h\bigr)}\,\times\ldots\\
 \hspace{-1.5cm}
  \ldots\times
  \frac{\prod\limits_{j_{N-m-1}=1\atop j_{N-m-1}\neq i_1}^{m+2}
  \bigl(\gamma_{m+1,\,i_1}-\gamma_{m+2,\,j_{N-m-1}}-\frac{\h}{2}\bigr)}
  {\prod\limits_{k_{N-m-1}=1\atop k_{N-m-1}\neq i_1}^m\!\!\!\!\!
  \bigl(\gamma_{m+1,\,i_1}-\gamma_{m+1,\,k_{N-m-1}}-\h\bigr)}\,\,
  \frac{\prod\limits_{j_{N-m}=1\atop j_{N-m}\neq i_1}^{m+1}
  \bigl(\gamma_{m,\,i_1}-\gamma_{m+1,\,j_{N-m}}-\frac{\h}{2}\bigr)}
  {\prod\limits_{k_{N-m}=1\atop k_{N-m}\neq i_1}^m\!\!\!\!\!
  \bigl(\gamma_{m,\,i_1}-\gamma_{m,\,k_{N-m}}\bigr)}\\
  \times
  \frac{\prod\limits_{j_{N+1-m}=1\atop j_{N+1-m}\neq i_1}^m
  \bigl(\gamma_{m-1,\,i_1}-\gamma_{m,\,j_{N-m}}-\frac{\h}{2}\bigr)}
  {\prod\limits_{k_{N+1-m}=1\atop k_{N+1-m}\neq i_1}^{m-1}\!\!\!\!\!
  \bigl(\gamma_{m-1,\,i_1}-\gamma_{m-1,\,k_{N+1-m}}-\h\bigr)}\,
  \,e^{-\h\sum\limits_{a=m-1}^{N-1}\pr_{a,\,i_1}}\Big\}\cdot\psi_R^{(m)}\,
  =\ldots
 \ee
Then we make cancelations due to the factors
$\prod\limits_{n=m+1}^{N-1}\prod\limits_{i,j\atop i\neq
j}^n\Gamma_1\Big(\gamma_{ni}-\gamma_{nj}\Big|\h\Big)$ in
\eqref{RightWhitVec} and relation \eqref{DeltaMeasureRelation}, take
into account the factor
$\prod\limits_{a=1}^{m-1}\prod\limits_{b=1}^m
  \Gamma_1\Big(\gamma_{m-1,\,a}-\gamma_{mb}+\frac{\h}{2}\Big|\h\Big)$
in \eqref{RightWhitVec} satisfying the relation:
 \be\label{ShiftGammaRelation3}
  e^{-\h\bigl(\pr_{m-1,\,j}+\pr_{m,\,i}\bigr)}\cdot
  \prod_{a=1}^{m-1}\prod_{b=1}^m\Gamma_1\Big(\gamma_{m-1,\,a}
  -\gamma_{n,\,b}+\frac{\h}{2}\Big|\,\h\Big)\\
  =\,\frac{\prod\limits_{r=1\atop r\neq j}^{m-1}\!\bigl(\gamma_{m-1,\,r}
  -\gamma_{m,\,i}-\frac{\h}{2}\bigr)}
  {\prod\limits_{p=1\atop p\neq i}^m\!\bigl(\gamma_{m-1,\,j}
  -\gamma_{m,\,p}-\frac{\h}{2}\bigr)}\,\,\,
  \prod_{a=1}^{m-1}\prod_{b=1}^m\Gamma_1\Big(\gamma_{m-1,\,a}
  -\gamma_{m,\,b}+\frac{\h}{2}\Big|\,\h\Big)\cdot
  e^{-\h\bigl(\pr_{m-1,\,j}+\pr_{m,\,i}\bigr)}\,,
 \ee
and then arrive to the following:
 \be
 \hspace{-1cm}
  \ldots=\,-\frac{1}{\h}\sum_{i_1=1}^m\,
  \prod\limits_{k_{N-m}=1\atop k_{N-m}\neq i_1}^m
  \frac{1}{\gamma_{m,\,i_1}-\gamma_{m,\,k_{N-m}}}\\
 \hspace{-1cm}
  \times\sum_{i_{N+1-m}=1}^{m-1}\!\!
  \frac{\prod\limits_{j=1\atop j\neq i_{N+1-m}}^{m-1}\!\!
  \bigl(\gamma_{m-1,\,j}-\gamma_{m,\,i_1}-\frac{\h}{2}\bigr)}
  {\prod\limits_{k_{N+1-m}=1\atop k_{N+1-m}\neq i_{N=1-m}}^{m-1}
  \bigl(\gamma_{m-1,\,i_{N+1-m}}-\gamma_{m-1,\,k_{N+1-m}}-\h\bigr)}\,\,
  e^{-\h\sum\limits_{a=m}^{N-1}\pr_{a,\,i_1}\,-\,\pr_{m-1,\,i_{N=1-m}}}
  \cdot\psi_R^{(m)}
 \ee
 \be
  =\,-\frac{1}{\h}\sum_{i_{N+1-m}=1}^{m-1}
  \prod\limits_{k_{N+1-m}=1\atop k_{N+1-m}\neq i_{N+1-m}}^{m-1}
  \frac{1}{\gamma_{m-1,\,i_{N+1-m}}-\gamma_{m-1,\,k_{N+1-m}}-\h}\\
  \times\,\underbrace{\sum_{i_1=1}^m\,
  \frac{\prod\limits_{j=1\atop j\neq i_{N=1-m}}^{m-1}\!\!
  \bigl(\gamma_{m-1,\,j}-\gamma_{m,\,i_1}-\frac{\h}{2}\bigr)}
  {\prod\limits_{k_{N-m}=1\atop k_{N-m}\neq i_1}^m\bigl(\gamma_{m,\,i_1}
  -\gamma_{m,\,k_{N-m}}\bigr)}}_{=0}\,\,
  e^{-\h\sum\limits_{a=m}^{N-1}\pr_{a,\,i_1}\,-\,\pr_{m-1,\,i_{N=1-m}}}
  \cdot\psi_R^{(m)}\,=\,0\,,
 \ee
where the sum above vanishes due to \eqref{Combin1}. This completes
the proof of Proposition.

\section{Mellin-Barnes integral and its asymptotic}

Now we are ready to construct the integral representation for the
$\Gr_{m,N}$-Whittaker function \eqref{GrWhittaker}. Let us
substitute \eqref{LeftWhitVec} and \eqref{RightWhitVec} into
\eqref{GrWhittaker} and then obtain:
 \be
  \Psi^{(m,\,N)}_{\underline{\gamma}_N}(x)\,\,
  =\,\,e^{-x\frac{m(N-m)}{2}}\int
  \prod_{n=1}^{N-1}d\underline{\gamma}_n\,\,
  \prod_{n=2}^{N-1}\prod_{i,j=1\atop i\neq j}^n
  \frac{1}{\Gamma_1\bigl(\gamma_{ni}-\gamma_{nj}\bigl|\h\bigr)}\,\,
  e^{-\frac{x}{\h}\sum\limits_{k=1}^m\gamma_{mk}}\\
  \times\delta(\gamma_{11})\prod_{i=1}^m
  \prod_{j=1}^N\Gamma_1\Big(\gamma_{N-1,\,i}-\gamma_{Nj}+\frac{\h}{2}
  \Big|\h\Big)\cdot
  \prod_{a=1}^{m-1}\prod_{b=1}^m
  \Gamma_1\Big(\gamma_{m-1,\,a}-\gamma_{mb}+\frac{\h}{2}\Big|\h\Big)\\
  \times\prod_{n=2\atop n\neq m}^{N-1}
  \Big[\delta\Big(\sum_{j=1}^n\gamma_{nj}
  -\sum_{i=1}^{n-1}\gamma_{n-1,\,i}\Big)
  \prod_{k=1}^{n-1}
  \delta\bigl(\gamma_{n-1,\,k}-\gamma_{nk}+\frac{\h}{2}\bigr)
  \prod_{i,j=1\atop i\neq j}^n\Gamma_1\bigl(\gamma_{ni}-\gamma_{nj}
  \bigl|\h\bigr)\Big]\\
  \times\,e^{\i\pi\ga_{11}}\prod_{i=1}^{m-1}\prod_{j=1}^m
  \frac{1}{\Gamma_1\bigl(\gamma_{m-1,\,i}
  -\gamma_{m,j}+\frac{\h}{2}\bigl|\h\bigr)}\\
  =\,e^{-x\frac{m(N-m)}{2}}\int\limits_{\CS}
  \prod_{n=1}^{N-1}d\underline{\gamma}_n\,\,
  e^{-\frac{x}{\h}\sum\limits_{k=1}^m\gamma_{mk}}\,\,\,
  \frac{\prod\limits_{i=1}^m\prod\limits_{j=1}^N
  \Gamma_1\Big(\gamma_{N-1,\,i}-\gamma_{Nj}+\frac{\h}{2}
  \Big|\h\Big)}
  {\prod\limits_{i,j=1\atop i\neq j}^m
  \Gamma_1\bigl(\gamma_{ni}-\gamma_{nj}\bigl|\h\bigr)}\\
  \times\,\delta(\gamma_{11})
  e^{\i\pi\ga_{11}}\prod_{n=2\atop n\neq m}^{N-1}
  \Big[\delta\Big(\sum_{j=1}^n\gamma_{nj}
  -\sum_{i=1}^{n-1}\gamma_{n-1,\,i}\Big)
  \prod_{k=1}^{n-1}
  \delta\bigl(\gamma_{n-1,\,k}-\gamma_{nk}+\frac{\h}{2}\bigr)\Big]\,
  =\ldots
 \ee
Making integration over
$\prod\limits_{n=1}^{N-2}d\underline{\gamma}_n
\prod\limits_{k=m+1}^Nd\gamma_{N-1,\,k}$ we integrate out the
delta-functions and arrive to
 \be\label{GrWhittakerCalculation}
 \hspace{-1cm}
  \ldots=\,e^{-x\frac{m(N-m)}{2}}\int
  \prod_{n=1}^{N-1}d\underline{\gamma}_n\,\,
  e^{-\frac{x}{\h}\sum\limits_{k=1}^m
  \bigl(\gamma_{N-1,k}-\frac{N-1-m}{2}\h\bigr)}\,\,\,
  \frac{\prod\limits_{i=1}^m\prod\limits_{j=1}^N
  \Gamma_1\Big(\gamma_{N-1,\,i}-\gamma_{Nj}+\frac{\h}{2}
  \Big|\h\Big)}
  {\prod\limits_{i,j=1\atop i\neq j}^m
  \Gamma_1\bigl(\gamma_{ni}-\gamma_{nj}\bigl|\h\bigr)}\,.
 \ee
Finally, we shift the integration contour by
$\ga_k=\ga_{N-1,\,k}+\frac{\h}{2}\,,k=1,\ldots,m\,,$ and readily get
\eqref{MBIntegralRep}.\, $\Box$

\subsection{Asymptotic of $\Psi^{(m,N)}(x,0,\ldots,0)$}

To complete the analysis of integral representation
\eqref{MBIntegralRep} let us derive its asymptotic when
$x\to-\infty$.

\emph{Proof of Theorem 2.2.} The contour of integration $\CC=\CC_m$
is a product of $m$ copies of contour $\CC_1$, corresponding to
integration over $\ga_k,\,k=1,\ldots,m$, going from $\e-\i\infty$ to
$\e+\i\infty$ with $\e>\max\{\la_1,\ldots,\la_N\}$. Let us enclose
$\CC_1$ by a half-circle of infinitely large radius in the left
half-plane, then the closed contour $\tilde{\CC}_1$ embraces all the
poles of gamma-factors
$$
 \prod_{j=1}^N\Gamma_1\bigl(\gamma_k-\la_j\bigr|\h\bigr)
$$
in \eqref{MBIntegralRep}. Actually, we can replace the integration
over $\CC_1$ by integration over $\tilde{\CC}_1$, since the
contribution over the half-circle of infinitely large radius
vanishes due to the exponent $e^{-\h^{-1}x\gamma_k}$ in the
integrand of \eqref{MBIntegralRep}. Then this transformation of
integration contour allows to calculate the integral as the sum over
the residues at poles of Gamma-factors. Namely, each Gamma-function
$\Gamma_1(\ga_k-\la_j|\h)$ has the poles at
$\ga_k=\la_j-n\h,\,n=0,1,2,\ldots$, and therefore integration over
$\ga_k$ implies the following:
$$
 e^{-x\ga_k}\prod_{j=1}^N\Gamma_1\bigl(\ga_k-\la_j\bigr|\h\bigr)\,
 =\,\sum_{i=1}^N\sum_{n=0}^{\infty}n!\,\,e^{-x(\la_j-n\h)}
 \prod_{j=1\atop j\neq i}^N
 \Gamma_1\bigl(\la_i-\la_k-n\h\bigr|\h\bigr)\,.
$$
When $x\to-\infty$ only the terms with $n=0$ give contributions into
the asymptotic:
$$
 e^{-x\ga_k}\prod_{j=1}^N\Gamma_1\bigl(\ga_k-\la_j\bigr|\h\bigr)\,\,
 \sim\,\,\sum_{i=1}^Ne^{-x\la_j}
 \prod_{j=1\atop j\neq i}^N
 \Gamma_1\bigl(\la_i-\la_j\bigr|\h\bigr)\,,
\hspace{1.5cm}
 x\to-\infty\,.
$$
At the next step we obtain:
 \be
  e^{-x(\ga_1+\ga_2)}\prod_{j=1}^N
  \Gamma_1\bigl(\ga_1-\la_j\bigr|\h\bigr)
  \Gamma_1\bigl(\ga_2-\la_j\bigr|\h\bigr)\\
  \sim\,\,\sum_{i_1=1}^N
  \prod_{j=1\atop j\neq i_1}^N
  \Gamma_1\bigl(\la_{i_1}-\la_j\bigr|\h\bigr)
  \sum_{i_2=1\atop i_2\neq i_1}^Ne^{-x(\la_{i_1}+\la_{i_2})}
  \frac{\prod\limits_{j=1\atop j\neq i_2}^N
  \Gamma_1\bigl(\la_{i_2}-\la_j\bigr|\h\bigr)}
  {\Gamma_1\bigl(\la_{i_1}-\la_{i_2}\bigr|\h\bigr)
  \Gamma_1\bigl(\la_{i_2}-\la_{i_1}\bigr|\h\bigr)}\\
  =\,\sum_{i_1=1}^N
  \prod_{j=1\atop j\neq i_1}^N
  \Gamma_1\bigl(\la_{i_1}-\la_j\bigr|\h\bigr)
  \sum_{i_2=1\atop i_2\neq i_1}^Ne^{-x(\la_{i_1}+\la_{i_2})}
  \frac{\prod\limits_{j=1\atop j\neq i_1,\,i_2}^N
  \Gamma_1\bigl(\la_{i_2}-\la_j\bigr|\h\bigr)}
  {\Gamma_1\bigl(\la_{i_1}-\la_{i_2}\bigr|\h\bigr)}\,.
 \ee
In this way we proceed step by step over $k$, making cancelations of
Gamma-factors in measure $\mu_m(\ug)$, and finally we arrive to
$\frac{N!}{m!(N-m)!}$ terms (with multiplicities $m!$) that can be
arranged into the $\frak{S}_N/\frak{W}_m$-orbit of the term
$$
 m!\,e^{-x(\ga_1+\ldots+\ga_m)}
 \prod_{i=1}^m\prod_{k=1}^{N-m}
 \Gamma_1\bigl(\la_i-\la_{m+k}\bigr|\h\bigr)\,.
$$
Thus we obtain \eqref{Asymptotics}. $\Box$

\vspace{1cm}

\noindent
{\sc Institute for Theoretical and Experimental Physics,\\
Bol. Cheremushkinskaya 25, Moscow 117218,\\
\emph{E-mail address:}\quad \verb"Sergey.Oblezin@itep.ru"}

\end{document}